\documentclass[a4paper,10pt]{article} 
\usepackage[applemac]{inputenc}
\usepackage{amssymb, amsmath, amscd}
\usepackage[pdftex]{graphicx}
\usepackage{charter}
\usepackage{mesmacrossansmathabxsansbabel}

\title{\textbf{Profinite completion and double-dual : isomorphisms and counter-examples}}
\date{}
\author{Colas Bardavid\footnote{\texttt{colas.bardavid@univ-rennes1.fr}}\\
IRMAR --- UMR 6625 du CNRS\\
Universit\'e de Rennes 1 Campus de Beaulieu \\
35042 Rennes CEDEX FRANCE }

\begin{document}
\maketitle
\bigskip

 \begin{center}  
\rule{13 cm}{.5pt} \end{center}
\begin{small}
\begin{center}
\begin{minipage}{13cm}
\textbf{Abstract --} We define, for any group $G$, finite approximations ; with this tool, we give a new presentation of the profinite completion $\widehat{\pi} : G \to \widehat{G}$ of an abtract group $G$. We then prove the following theorem : if $k$ is a finite prime field and if $V$ is a $k$-vector space, then, there is a natural isomorphism between $\widehat{V}$ (for the underlying additive group structure) and the additive group of the double-dual $V^{**}$. This theorem gives counter-examples concerning the iterated profinite completions of a group. These phenomena don't occur in the topological case.
\end{minipage}
\end{center}
\end{small}
 \begin{center}  
\rule{13 cm}{.5pt} \end{center}

\sct{Introduction. }
In this paper, we study the profinite completion of a certain class of groups, namely, the additive groups of vector spaces over $\F{p}$. The principal result is that, in this case, the profinite completion equals the double-dual. This study is based on a ‘‘dual’’ definition of the profinite completion of a group.

\pause

\sct{Brief survey of the classical point of view for profinite completion. }
As explained in \cite{galoiscoh} or \cite{ribeszalesskii}, one usually defines the profinite completion of a group\footnote{If we set in thecategory of \emph{topological} groups, we should precise : ‘‘... of a \emph{discrete} group $G$’’.} $G$ as follows. The profinite completion $\widehat{G}$ of $G$ is the projective limit (\ie the inverse limit) of the finite quotients of $G$ :
$$
\widehat{G}= \mathop{\varprojlim}\limits_{\begin{subarray}{c}
N\lhd G \\ [G:N]<\infty 
\end{subarray} }G/N.
$$

There is a more explicit form for this definition. Indeed, if $N, M$ are two normal subgroups of $G$ with $N\subset M$, we have a natural factorisation $\varphi_{N\subset M}$ of the canonical projection $\pi_{M}$ :
$$
\xymatrix@R=2mm{
& G/N\ar[dd]^-{\varphi_{N\subset M}}\\
G \ar[ur]^-{\pi_{N}} \ar[rd]_-{\pi_{M}} \\
&G/M
}
$$
One  can then write :
$$\widehat{G} = \left \{ (x_{N}) \in \prod_{\begin{subarray}{c}
N\lhd G \\ [G:N]<\infty 
\end{subarray}} G/N \tq \fa N\subset M,\quad \varphi_{N\subset M} \left ( x_{N} \right )= x_{M} \right \}.$$

\pause

\sct{Finite approximations and profinite completion. }
In this paper, we will use a ‘‘dual’’ (but equivalent) point of view for the profinite completion of a group. To begin with, we introduce the notion of ‘‘finite approximation’’, which will lead naturally to the concept of profinite completion.

\begin{itshape}
\ssct{Definition.} If $G$ is a group, we call \emph{finite approximation of $G$} every couple $\mathbf{v}=(F, \varphi)$ where $F$ is a finite group and $\varphi : G\to F$ a morphism. We denote $F=F_{\mathbf{v}}$ and $\varphi=\varphi_{\mathbf{v}}$. We say that $f : \mathbf{v}\to \mathbf{v'}$ is a morphism between $\mathbf{v}$ and $\mathbf{v'}$ if it is an arrow that makes the following diagram commute :
$$
\xymatrix@R=2mm{
&F_{\mathbf{v}}\ar[dd]^{f}\\
G\ar[ur]^{\varphi_{\mathbf{v}}}\ar[dr]_{\varphi_{\mathbf{v'}}}\\
&F_{\mathbf{v'}}}.
$$
We denote $\textbf{App}_{f}(G)$ the category of finite approximations of $G$.
\end{itshape}
\bigskip

Intuitively, a finite approximation of $G$ allows the mathematician to get some information about $G$ by only dealing with finite objects. Here are some examples, from various aeras of mathematics, of finite approximations :

\bigskip
\begin{itemize}
\item[\textbf{a) }] $\fonction{\R^*}{\Z/2 \Z}{x}{\text{sgn}(x)}$ the sign of a real number.
\bigskip
\item[\textbf{b) }] The reduction modulo $n$, $\Z\to \Z/n\Z$ and all the derived morphisms and generalizations, such that $\Z_{(p)}\to \Z/p \Z$, such that $GL_{m}(\Z)\to GL_{m}\left ( \Z/n \Z \right )$ or such that $\mathcal{O}_{K}\to \mathcal{O}_{K}/\mathfrak{P}$ if $K$ is a number field ;
\bigskip
\item[\textbf{c) }] If $X$ a topological space with a finite number of connected components, we can consider the ‘‘trace’’ on $\pi_{0}(X)$ of an automorphism : $\fonction{\Aut (X)}{\mathfrak{S}_{\pi_{0}(X)}}{\phi}{\pi_{0}(\phi)}$.
\bigskip
\item[\textbf{d) }] If we denote $\mathfrak{S}_{(\N)}=\limind_{n}\mathfrak{S}_{n}$ the group of permutation of $\N$ with finite support, we can still define a signature $\mathfrak{S}_{(\N)}\to \Z/2\Z$.
\bigskip
\item[\textbf{e) }] Finally, if $K/\Q$ is a Galois extension, then $\fonction{\Gal\!\!\left(\barre{\Q}/\Q \right )}{\Gal\!\!\left(K/\Q \right )}{\sigma}{\rst{\sigma}{K}}$ is a finite approximation.
\end{itemize}

\bigskip
\ssct{Profinite completion.} Then, one can define very naturally the profinite completion of $G$ as the projective limit of all the finite approximations of $G$. More precisely, (and without dealing with any problem of set theory)
$$
\widehat{G}=\left \{ \left ( g_{\mathbf{v}}\right )_{\mathbf{v}\in \mathbf{App}_f(G)}\in \prod_{\mathbf{v}}F_{\mathbf{v}} \tq \fa \psi : \mathbf{v}\to\mathbf{w},\,\, \psi(g_{\mathbf{v}})=g_{\mathbf{w}}\right \}
$$
which comes with the \emph{profinite projection}
$$
\widehat{\pi} : \fonction{G}{\widehat{G}}{g}{\left ( \varphi_{\mathbf{v}}(g)\right )_{\mathbf{v} \in  \mathbf{App}_f(G)}}.
$$
Intuitively, this object is what remains from $G$ when one can only deal with information of finite type ; some elements will be identified but, at the same time, some new elements will appear. Formally, in general, $\widehat{\pi}$ is not surjective or injective.

\ssct{Surjective finite approximations.} Among the finite approximations, some are surjective ; they form a full subcategory $\mathbf{App}_{f}^{s}(G)$ of $\mathbf{App}_{f}(G)$. In the same way that we have defined the profinite completion, we can then define the ‘‘surjective’’ profinite completion $$\widehat{G}^{s}=\mathop{\varprojlim}\limits_{\mathbf{v}\in \mathbf{App}_{f}^{s}(G)} F_{\mathbf{v}}.$$
The important fact about this object is that we have the following fact, whose proof is not difficult.

\begin{itshape}
\ssct{Proposition.}The natural morphism $\widehat{G}\to \widehat{G}^{s}$ is an isomorphism.
\end{itshape}

\pause

\ssctn
\sct{Profinite completion of the additive group of a vector space over $\F{p}$.}\label{compa}

\ssct{Profinite completion and double-dual.} Before looking at what happens in the situation where the base field is $\F{p}$, let us remark that, in the general case, there is a morphism of comparison between the profinite completion of an ‘‘additive group’’ and its double-dual. Let $k$ be a finite field and $V$ a vector space over $k$. We still denote by $V$ the underlying additive group.

\smallskip
Let $f$ and $g$ be two linear forms of $V$ and let $\lambda\in k$. The forms $f$, $g$ and $f+\lambda \cdot g$ are, in particular, finite approximations of $V$ (in the additive group of $k$) and we denote by $\mathbf{v}_{f}$, $\mathbf{v}_{g}$ and $\mathbf{v}_{f+\lambda \cdot g}$ the corresponding approximations. Now, let $\mathbf{x}=\left ( x_{\mathbf{v}} \right )_{\mathbf{v}}\in \widehat{V}$ be a ‘‘profinite’’ element.

\begin{itshape}
\ssct{Fact.} $x_{\mathbf{v}_{f+\lambda \cdot g}}=x_{\mathbf{v}_{f}}+\lambda \cdot x_{\mathbf{v}_{g}}$.
\end{itshape}

\smallskip
\proof{Indeed, we have the following diagram of morphisms of finite approximations
$$
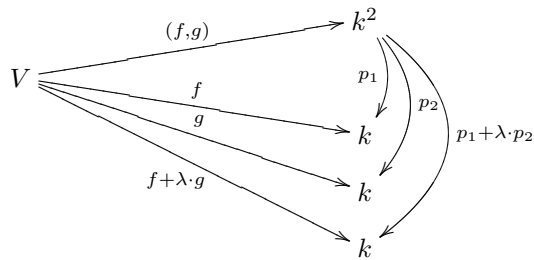
\xymatrix@C=4cm@R=3mm{
&k^2\ar@/^3mm/[dd]_{p_{1}}  \ar@/^6mm/[ddd]^{p_{2}}  \ar@/^11mm/[dddd]^{p_{1}+\lambda \cdot p_{2}}\\
V\ar[ur]^{(f,g)}\ar[dr]^-{f}\ar[ddr]^-{g}\ar[dddr]_-{f+\lambda \cdot g}\\
&k\\
&k\\
&k
}.$$
Then, if we denote by $\mathbf{w}$ the approximation $\fleche{(f,g)}{V}{k^2}$, the definition of the profinite completion imposes that
$
x_{\mathbf{v}_{f}}=p_{1}\left (  x_{\mathbf{w}}\right )$ and $ x_{\mathbf{v}_{g}}=p_{2}\left (  x_{\mathbf{w}}\right )$  and $x_{\mathbf{v}_{f+\lambda \cdot g}}=\left ( p_{1}+\lambda \cdot p_{2} \right )\left (  x_{\mathbf{w}}\right )$, that is 
$$x_{\mathbf{v}_{f+\lambda \cdot g}}=x_{\mathbf{v}_{f}}+\lambda \cdot x_{\mathbf{v}_{g}}.$$
}

Using this fact, one can define the morphism of comparison :
$$
\Psi : \xymatrix@R=2mm{
\widehat{V}\ar[r]& V^{**}\\
(x_{\mathbf{v}})_{\mathbf{v}}\ar@{|->}[r]& \left ( \fonction{V^*}{k}{f}{x_{\mathbf{v}_{f}}} \right )
}
$$

\ssct{The case where $k=\F{p}$.}
From now on, $p$ is a prime number and $k=\F{p}$. The interesting case is when $V$ is of infinite dimension. A good way to understand what happens is to consider $V=\left ( \Z/2\Z \right )^{\N}$. 

\smallskip

The first thing to do is to see that if $\varphi : V\to F$ is a finite \emph{surjective} approximation, then $F$ is isomorphic to (the additive group of) ${\left ( \F{p} \right )}^n$ for some $n$. Indeed, first of all, since $F$ is the homomorphic image of $V$, $F$ is abelian. Moreover, all the elements of $F$ satisfy $x^p=e$. Thus, the classification of the abelian finite groups gives the conclusion. 

\bigskip
We can now prove :

\begin{itshape}
\ssct{Theorem}\label{theo1}
Let $V$ be a vector space over $\F{p}$. Then, $\Psi : \widehat{V}\to V^{**}$ is an isomorphism.
\end{itshape}
\bigskip

\proof{We first prove that $\Psi$ is injective : let $\mathbf{x}=(x_{\mathbf{v}})_{\mathbf{v}}\in \widehat{V}$ such that for all linear form $f : V\to k$, $x_{\mathbf{v}_{f}}=0$. Let $\mathbf{v}$ be a finite surjective approximation of $V$ ; we can suppose that $\mathbf{v}=(k^n, \varphi)$, where $\varphi : V \to k^n$ is any morphism. By composing $\varphi$ with the $n$ projections $p_{i}$ to the factors $k$, one obtain $n$ morphisms. If we prove that the $n$ corresponding elements are equal to $0$, then, it will follow that $x_{\mathbf{v}}$ is equal to $0$ and, thus, that $\Psi$ is injective. But, and it is the (easy) key point, a morphism $V\to k$ of groups is actually a linear form, since we can rewrite the condition $\varphi(\lambda\cdot \vec{v})=\lambda\cdot \varphi(\vec{v})$ as $\varphi(\vec{v}+\cdots +\vec{v})= \varphi(\vec{v})+ \cdots + \varphi(\vec{v})$, for our base field is $\F{p}$. And, by assumption, all the $x_{\mathbf{v}_{f}}=0$. 
\smallskip

For the surjectivity, let $\Theta\in V^{**}$ be a double-dual element. We would like to find a profinite element $\mathbf{x}=(x_{\mathbf{v}})_{\mathbf{v}}\in \widehat{V}$ such that, for all linear form $f$ of $V$, one have $x_{\mathbf{v}_{f}}=\Theta(f)$. So, let $\mathbf{v}=(k^n, \varphi)$ (as we can suppose it) be a finite approximation of $V$. Let denote $p_{1}, \cdots, p_{n}$ the $n$ projections of $k^n$ to the factors $k$. Naturally, we define $x_{\mathbf{v}}$ by reconstructing it from the linear forms $p_{i}\circ \varphi$ :
$$
x_{\mathbf{v}}:=\left ( \Theta \left ( p_{1}\circ \varphi \right ), \Theta \left ( p_{2}\circ \varphi \right ), \ldots, \Theta \left ( p_{n}\circ \varphi \right ) \right )\in k^n.
$$
Now, we just have to check that the family $(x_{\mathbf{v}})$ is ‘‘compatible``. So let $\mathbf{v}=(k^n, \varphi)$ and $\mathbf{w}=(k^n, \psi)$ be two finite approximations and $g$ a morphism between them :
$$
\xymatrix@R=2mm{
&k^n\ar[dd]^{g}\\
V\ar[ur]^{\varphi}\ar[dr]_{\psi}\\
&k^m}.
$$
We want to prove that $g(x_{\mathbf{v}})=x_{\mathbf{w}}$. By composing with the $m$ projections $q_{i}$ of $k^m$, it suffices to prove it in the case where $m=1$ :
$$
\xymatrix@R=2mm{
&k^n\ar[dd]^{g}\\
V\ar[ur]^{\varphi}\ar[dr]_{\psi}\\
&k^m\ar[r]^{q_{1}}\ar[ddr]_{q_{m}}&k\\
&&\vdots\\
&&k}.
$$
So, we are brought to this situation
$$
\xymatrix@R=2mm{
&k^n\ar[dd]^{g}\\
V\ar[ur]^{\varphi}\ar[dr]_{\psi}\\
&k},
$$
where we know that $g$ can be written as $g=\lambda_{1}\cdot p_{1}+\cdots +\lambda_{n}\cdot p_{n}$, with $\lambda_{i}\in k$. The fact that the previous diagram commutes tells us that $\psi=\sum_{i} \lambda_{i} \cdot \left ( p_{i}\circ\varphi \right )$ ; and, now :
\begin{eqnarray*}
g(x_{\mathbf{v}})&=&g\left ( \left ( \Theta \left ( p_{1}\circ \varphi \right ), \Theta \left ( p_{2}\circ \varphi \right ), \ldots, \Theta \left ( p_{n}\circ \varphi \right ) \right ) \right )\\
&=& \sum_{i} \lambda_{i} \cdot \Theta \left ( p_{i}\circ \varphi \right )\\
&=& \Theta \left (  \sum_{i} \lambda_{i} \cdot \left ( p_{i}\circ\varphi \right )\right ) = \Theta (\psi)\\
&=& x_{\mathbf{w}},
\end{eqnarray*}
which concludes the proof. }
\bigskip

\ssct{$\widehat{\pi}$ and the canonical injection $i : V\to V^{**}$.} We denote $i : V\to V^{**}$ the canonical injection defined by $i(\vec{v})(f)=f(\vec{v})$. One can improve a bit the theorem \textbf{\ref{compa}.\ref{theo1}} :  the isomorphism $\Psi$ between $\widehat{V}$ and $V^{**}$ through $\Psi$ identifies $\widehat{\pi}$ with $i$. The proof is easy.

\bigskip

\begin{itshape}
\ssct{Theorem.}
Let $V$ be a vector space over $\F{p}$. Then, $\Psi : \widehat{V}\to V^{**}$ is an isomorphism and the diagram
$$
\xymatrix@R=2mm{
&\widehat{V}\ar[dd]^{\Psi}_{\wr}\\
V\ar[ru]^{\widehat{\pi}} \ar[dr]_{i}\\
& V^{**}
}
$$
commutes.
\end{itshape}

\bigskip
\ssct{Remark. }\label{referee}One can prove the theorem \textbf{\ref{compa}.\ref{theo1}} with more abstracted arguments. To begin with, we know (\cf for example \cite[\S 2.7, theorem 2]{livrecategories}) that, in a general category $\mathcal{C}$, if the limits exist, we always have the natural isomorphism
$$
\Hom \left ( \limind _{i} X_{i}, X \right )  \cong \limproj_{i} \Hom \left (X_{i}, X \right ) .
$$
Moreover, in the case of $k$-vector spaces, this isomorphism is linear ; thus, for a system of $k$-vector space $V_{i}$, we have :
$$
\left (  \limind_{i} V_{i} \right )^* \cong \limproj_{i} \left ( {V_{i}}^* \right ).
$$
Let $k$, from now on, be a field and $V$ a $k$-vector space. If we denote by $\left ( Y_{i} \right )_{i}$ the system of finite-dimensional subvector spaces of $V^*$, we have $V^*=\limind_{i} Y_{i}$ and thus, thanks the previous isomorphism :
$$
V^{**}\cong \limproj_{i} \left ( {Y_{i}}^* \right ).
$$
Moreover, there is a natural bijection between the finite-dimensional subspaces of $V^*$ and the finite-codimensional subspaces of $V$, \emph{via} the application 
$$Y\mapsto Y^\perp :=\left \{v\in V \tq \fa \varphi\in Y, \,\varphi(v)=0\right \}. $$
Besides, if $Y$ is a finite-dimensional subspace of $V^*$ then the dual $Y^*$ is naturally isomorphic to $V/{Y^{\perp}}$. Consequently, if we denote by $(Z_{j})_{j}$ the system of finite-codimensional subspaces of $V$, we have :
$$
V^{**}\cong \limproj_{j}\left ( V/Z_{j} \right ).
$$
But, if $k=\F{p}$ for a prime number $p$, one can identify the $k$-vector space $V$ with its underlying additive group\footnote{We denote $\omega : k-\text{\textbf{Vs}}\to \Gr$ the forgetful functor from the category of $k$-vector spaces to the category of groups.} $\omega(V)$, its dual $V^*$ with $\Hom_{\Gr}\left ( \omega(V), \omega(\F{p}) \right )$, and its finite dimensional quotients with the finite quotient of  $\omega(V)$. We thus finally get the expected alternative proof of the theorem \textbf{\ref{compa}.\ref{theo1}}.

\pause

\ssctn
\sct{A family of counter-examples.}
One would like to know if, given a group $G$, one have $\widehat{\widehat{G}}\simeq \widehat{G}$. This fact is known to be false (\cf example 4.2.13 of \cite{ribeszalesskii}), but as we will see, it is still false, in general, after taking $i$ times the profinite completion.

\bigskip
\ssct{The sequence of $i$-th profinite completions. }We introduce the following notation. If $G$ is a group, we denote 
 ${\widehat{G}}^{[1]}= \widehat{G}$
and  ${\widehat{G}}^{[i+1]}= \widehat{  
{\widehat{G}}^{[i]}
   }$. These groups come with projections, as follows :
   $$
   \xymatrix{
   G\ar^{{\widehat{\pi}}^{[1]}}[r] & {\widehat{G}}^{[1]} \ar^{{\widehat{\pi}}^{[2]}}[r]& {\widehat{G}}^{[2]}\ar[r]&\cdots \ar[r]& {\widehat{G}}^{[i]}\ar^{{\widehat{\pi}}^{[i+1]}}[r] & {\widehat{G}}^{[i+1]} \ar[r]&\cdots
   }.
   $$
We will prove that, in general, none of the ${\widehat{\pi}}^{[i]}$ is an isomorphism.

\bigskip

\begin{itshape}
\ssct{Proposition.}\label{contrex}Let $p$ be a prime number and $k=\F{p}$. Let $V$ be (the additive group of) a $k$-vector space of infinite dimension. Then, in the following sequence
$$
   \xymatrix{
   V\ar^{{\widehat{\pi}}^{[1]}}[r] & {\widehat{V}}^{[1]} \ar^{{\widehat{\pi}}^{[2]}}[r]& {\widehat{V}}^{[2]}\ar[r]&\cdots \ar[r]& {\widehat{V}}^{[i]}\ar^{{\widehat{\pi}}^{[i+1]}}[r] & {\widehat{V}}^{[i+1]} \ar[r]&\cdots
   }
   $$
   all the ${\widehat{\pi}}^{[i]}$ are injective but non-surjective morphisms.
\end{itshape}

\bigskip

\proof{
This follows from the identification of the arrows ${\widehat{\pi}}^{[i]}$ with the canonical injections of a vector space in its double-dual, and from the fact that these injections are injective but non-surjective when the vector spaces are of infinite dimension, cf. Th\'eor\`eme 6, \S 7, n$^\circ$5 of \cite{bbkialglin}.
}

\pause
\newpage
\ssctn
\sct{Conclusion : abstract setting vs. topological setting.}\label{conclu}
This study has been given for groups but a similar point of view can  be applied to \emph{topological} groups. In this case, we start with a topological group $\mathcal{G}$ and we consider the category $\textbf{App}_{discr}(G)$ of finite and \emph{discrete} approximations : they are couples $\mathbf{v}=(F, \varphi)$, where $F$ is a discrete and finite topological group and $\varphi : \mathcal{G}\to F$ a continuous morphism of groups.

\bigskip
One obtain the (topological) profinite completion of $\mathcal{G}$, wich is, as well-known, a topological group, compact and totally disconnected (cf. \cite{galoiscoh}), and one obtain a profinite projection, which is a continuous morphism :
$$
\widehat{\pi}^{top} : \mathcal{G}\to \widehat{\mathcal{G}}^{top}.
$$
More generally, as previously done, one can define the sequence of iterated (topological) profinite completions : 
 $$
   \xymatrix@C=1.2cm{
   \mathcal{G}\ar^{{\widehat{\pi}}^{[1], top}}[r] & {\widehat{\mathcal{G}}}^{[1], top} \ar[r]&\cdots \ar[r]& {\widehat{\mathcal{G}}}^{[i], top}\ar^{{\widehat{\pi}}^{[i+1], top}}[r] & {\widehat{\mathcal{G}}}^{[i+1], top} \ar[r]&\cdots
   }.
   $$
The situation is then totally different than before. Indeed, we have :

\bigskip

\begin{itshape}
\ssct{Proposition. }\label{isotopo}Let $\mathcal{G}$ be a topological group. Then, for all $i\geq 2$, the arrows ${{\widehat{\pi}}^{[i], top}}$ are isomorphisms of topological groups.
\end{itshape}

\bigskip

\ssct{Profinite groups : abstract setting and topological setting.} There is a synthetical way to see the fundamental difference between the propositions \textbf{\ref{compa}.\ref{contrex}} and \textbf{\ref{conclu}.\ref{isotopo}}. For this sake, we introduce two notions of profinite groups. We will say that a group $G$ is \emph{profinite} if it is the projective limit of a system of finite groups ; we will say that a topological group $\mathcal{G}$ is \emph{topologically profinite} if it is the projective limit of a system of finite and discrete groups. We then have :

\bigskip

\begin{itshape}
\ssct{Theorem} Let $\mathcal{G}$ be a topological group. Then :
$$ \mathcal{G}\text{ is topologically profinite}\ssi \widehat{\pi}^{top} : \mathcal{G}\to \widehat{\mathcal{G}}^{top}\text{ is an isomorphism}.$$
\end{itshape}

\begin{itshape}
\ssct{Proposition} Let $G$ be a group. Then :
$${G}\text{ is profinite}\,\Leftarrow\, \widehat{\pi} :{G}\to \widehat{{G}}\text{ is an isomorphism}$$
$${G}\text{ is profinite}\,\nRightarrow\, \widehat{\pi} :{G}\to \widehat{{G}}\text{ is an isomorphism}.$$
\end{itshape}

\ssct{A positive answer.} One could legitimately be disapointed by the non-equivalence of $G$ being profinite and of $\widehat{\pi}$ being an isomorphism. Indeed, on the one hand, there is the very classical definition of a profinite group and, on the other hand, there is the deep property for a group to have its profinite projection $\widehat{\pi}$ to be an isomorphism (such a group, in a way, is separated ---~for $\widehat{\pi}$ is injective~--- and complete ---~for $\widehat{\pi}$ is surjective). One would have expected these two to coincide...

\bigskip
Fortunately, there is a positive result in this direction. It is a difficult result, which has been published in 2007 by Nikolay Nikolov and Dan Segal, \cf \cite{niknik1} and \cite{niknik2}, and whose proof uses the classification of finite simple groups. In order to state their result, let us remark that if $G$ is an (abstract) profinite group, if we write $G=\limproj_{i} F_{i}$, where the $F_{i}$'s are finite, and if we endow each of the $F_{i}$'s with the discrete topology, then we can view $G$ as a topological group.

\bigskip
\begin{itshape}
\ssct{Theorem}
Let $G$ be an (abstract) profinite group, which is topologically of finite type for the associated topology. Then, $\widehat{\pi}:{G}\to \widehat{{G}}$ is an isomorphism.
\end{itshape}

\ppause

\ssct{Acknowledgments. }My first acknowledgments go to Xavier Caruso for many helpful discussions. I would like also to thank the referee for many valuable comments and for making me known the alternative proof $\textbf{\ref{compa}.\ref{referee}}$.

\bibliography{bibliographie}

 \end{document}